 \documentclass[onefignum,onetabnum]{siamart171218}
\usepackage{subeqnarray}

\usepackage{lipsum}
\usepackage{amsfonts}
\usepackage{graphicx}
\usepackage{epstopdf}
\usepackage{algorithmic}
\ifpdf
  \DeclareGraphicsExtensions{.eps,.pdf,.png,.jpg}
\else
  \DeclareGraphicsExtensions{.eps}
\fi

\newsiamremark{remark}{Remark}
\newsiamremark{hypothesis}{Hypothesis}
\crefname{hypothesis}{Hypothesis}{Hypotheses}
\newsiamthm{claim}{Claim}

\headers{Singularities of hyperbolic PDEs}{J. Eshima, L. Deike, and H. A. Stone}

\title{Similarity solutions of shock formation for first-order strictly hyperbolic systems \thanks{Submitted to the editors DATE.
\funding{This research was conducted under the NSF Grant No. 2246791 to H.A.S. and NSF grant No. 2242512 to L.D., and the Eli and Britt Harari Fellowship in the Department of Mechanical and
Aerospace Engineering, Princeton University to J.E.}}}

\author{Jun Eshima\thanks{Department of Mechanical and Aerospace Engineering, Princeton University, Princeton, New Jersey 08544, USA (\email{jeshima@princeton.edu}, \email{ldeike@princeton.edu}, \email{hastone@princeton.edu}).}
\and Luc Deike \footnotemark[2]~\thanks{High Meadows Environmental Institute, Princeton University, Princeton, New Jersey 08544, USA.}
\and Howard A. Stone\footnotemark[2]}

\usepackage{amsopn}

\makeatletter
\newcommand*{\addFileDependency}[1]{
  \typeout{(#1)}
  \@addtofilelist{#1}
  \IfFileExists{#1}{}{\typeout{No file #1.}}
}
\makeatother

\newcommand*{\myexternaldocument}[1]{%
    \externaldocument{#1}%
    \addFileDependency{#1.tex}%
    \addFileDependency{#1.aux}%
}

\myexternaldocument{supplement}

\begin{document}

\maketitle

\begin{abstract}
  Shocks due to hyperbolic partial differential equations (PDEs) appear throughout mathematics and science. The canonical example is shock formation in the inviscid Burgers' equation $\partial u/\partial t + u \partial u/\partial x = 0$. Previous studies have shown that when shocks form for the inviscid Burgers' equation, for positions and times close to the shock singularity, the dynamics are locally self-similar and universal, i.e., dynamics are equivalent regardless of the initial conditions. In this paper, we show that, in fact, shock formation is self-similar and universal for general first-order strictly hyperbolic PDEs in one spatial dimension, and the self-similarity is like that of the inviscid Burgers' equation. An analytical formula is derived for the self-similar universal solution.
\end{abstract}

\begin{keywords}
  Hyperbolic PDEs, Similarity solutions, Burgers' equations
\end{keywords}

\begin{AMS}
  35L67, 35L60, 35C06, 76L05
\end{AMS}

\section{Introduction}\label{sec:intro}

Shocks due to hyperbolic partial differential equations (PDEs) have long been studied in mathematics and science with seminal works given by Challis \cite{challis48}, Stokes \cite{Stokes48}, and Riemann \cite{Riemann60}. For a summary, see \cite{Dafermos, JohnsonCheret}. The sustained interest arises due to the prevalence of such equations, with  applications spanning many fields of physics, such as fluid dynamics, including the shallow water and compressible gas equations, soft matter physics, such as nonlinear elasticity, magnetohydrodynamics, as occurs in astrophysics and plasma physics, as well as more generic  phenomena such as traffic flow.  

Consider the system of PDEs given by
\begin{equation}
    \frac{\partial \boldsymbol{f}}{\partial t}= \boldsymbol{M}(\boldsymbol{f}) \cdot  \frac{\partial \boldsymbol{f}}{\partial x},\label{eq:hyperbolic_PDE}
\end{equation}
where $\boldsymbol{f}(x,t)=(f_1(x,t), \cdots, f_N(x,t))$ are the $N$ dependent variables, $x$ is the spatial variable, $t$ denotes time, and $\boldsymbol{M}(\boldsymbol{f})$, an $N \times N$ matrix, is some smooth function of $\boldsymbol{f}$. A system of PDEs given by (\ref{eq:hyperbolic_PDE}) is said to be \textit{hyperbolic} \cite{Dafermos} if $\boldsymbol{M}$ has a basis of eigenvectors, i.e., is diagonalizable, with only real eigenvalues. Furthermore, if the eigenvalues are distinct, then the system is said to be \textit{strictly hyperbolic}. In this work, we consider strictly hyperbolic PDE systems with one space dimension. 

Mathematical studies of shocks often analyze fundamental questions such as uniqueness, existence, and well-posedness. A less commonly asked question is how the shocks are formed, which is the main focus of this paper. Understanding how shocks form is of practical value, as it naturally leads to the understanding of their subsequent regularization. To answer such a question, a novel approach was presented by Pomeau \cite{Pomeau08} and Eggers, and Fontelos \cite{Eggers09}, who investigated the inviscid Burgers' equation
\begin{equation}
    \frac{\partial u}{\partial t}+u\frac{\partial u}{\partial x} = 0, \label{eq:burger}
\end{equation}and showed that the shock formation is locally self-similar through a universal solution (see \cref{subsec:background}). The Burgers' equation is often regarded as the prototypical equation for shock formation for various applications \cite{Bonkile18}, such as benchmarking numerical codes. The simplicity of the Burgers' equation also leads to its appearance in unexpected areas of mathematics, including probabilistic models such as random matrix theory \cite{Menon12}.
From the perspective of shock formation, several authors \cite{Eggers17, Ryutovl19, Graziani22, Lai25, Eshima25} have shown that particular examples of shock formation in hyperbolic PDEs are closely related to the Burgers' equation. 

The formation of singularities from smooth conditions, such as those that arise from the Burgers' equation, is also of interest when investigating the well-posedness of flow equations. The most famous example is whether the Navier-Stokes equations form singularities in finite-time, which is one of the Millenium prize problems \cite{Fefferman06}. Notably, the identification of similarity solutions of finite-time singularities have been used in approaching such problems \cite{Eggers19, Wang23}.   

In this paper, we show that the relation to the Burgers' equation is a generic property of shock formation of hyperbolic PDEs. In particular, shock formation of first-order systems of strictly hyperbolic PDEs in one dimension is locally self-similar and we derive the self-similar solutions.  

To clarify, in our analysis, we assume that a finite-time shock exists and investigate the local self-similarity of such shock formation. We therefore do not investigate whether shocks form or not. For such studies, see for example \cite{Lax64,John74}. In particular, John \cite{John74} proved that a finite-time shock must form when the initial data for strictly hyperbolic PDEs (\ref{eq:hyperbolic_PDE}) is sufficiently small, i.e., the smallness leads the waves to behave as simple waves, which therefore form a Burgers-type shock.

The method shown in this paper is a generalization of the analysis of the Burgers' equation by Pomeau \cite{Pomeau08} and Eggers and Fontelos \cite{Eggers09}, where there is one dependent variable $u$ (\ref{eq:burger}). Additionally, it is a generalization of a result in a recent work by the authors of this paper \cite{Eshima25}, where  the same self-similar solutions of shock formation are derived for the specific example of thin film Marangoni flow, where there are three dependent variables. In this paper, we consider the broad class of strictly hyperbolic PDEs with $N$ dependent variables. 

The structure of the paper is as follows. In \cref{subsec:background} the existing self-similar shock solution for the Burgers' equation is given  and in \cref{subsec:main} we state the main, analogous result of this paper for equation (\ref{eq:hyperbolic_PDE}). The result is derived in detail in \cref{sec:deriv}, with some corollaries provided in \cref{sec:corollaries}. Finally, in \cref{sec:verif}, the method of solution is presented through an example using the shallow water equations, where the numerical solution provides independent confirmation of the theoretical prediction.

\section{Background: Burgers' Equation}\label{subsec:background}

The similarity analysis of the inviscid Burgers' equation forms the basis of the results to follow.  Thus, we first remind a reader of the similarity solution of the inviscid Burgers' Equation (\ref{eq:burger}), following the approach given in \cite{Eggers09}. 

Let a shock form at $(x_*,t_*)$ with $u(x_*,t_*) = u_*$.  Then, as the shock forms, for positions and times close to the shock, i.e., $|x-x_*| \ll 1$ and $t \rightarrow t_*^-$, the leading-order expansion of (\ref{eq:burger}) gives 
\begin{equation}
    \frac{\partial u}{\partial t}+u_*\frac{\partial u}{\partial x}=0,\label{eq:advection}
\end{equation}
which is a linear advection equation of velocity $u_*$ and does not develop a shock by itself. 

For further analysis, the next-order nonlinearity is considered, where the shock occurs. We also consider variables local to the shock formation. Let $\overline{x} :=x-x_* +u_*\tau$ be the frame of the shock and let $\tau:=t_*-t$ be the local time, where the sign is chosen so that $\tau$ remains positive as the shock is approached $t \rightarrow t_*^-$. Finally, we let $u':=u-u_*$. Then, (\ref{eq:burger}) becomes
\begin{equation}
    \frac{\partial u'}{\partial \tau}-u'\frac{\partial u'}{\partial \overline{x}}=0.\label{eq:burgers_shockframe}
\end{equation}

The scalings for the self-similar ansatz may be derived by balancing $\partial u'/\partial \tau$ and $u'\partial u'/\partial \overline{x}$ to deduce that $\overline{x} = \textit{O}(\tau^{\alpha})$ and $u' = \textit{O}(\tau^{\alpha-1})$ as $\tau \rightarrow 0^+$ for some constant $\alpha$, to be determined. Explicitly, the self-similar ansatz is then given by
\begin{equation}
    u' = \tau^{\alpha-1}F(\xi),\label{eq:ansatz_burgers}
\end{equation}
where $\xi :=\overline{x} \tau^{-\alpha}$ is the similarity coordinate. The ansatz (\ref{eq:ansatz_burgers}) may be substituted into (\ref{eq:burgers_shockframe}) to give an ordinary differential equation (ODE):
\begin{equation}
    (\alpha -1) F -\alpha \xi \frac{dF}{d\xi}-F\frac{dF}{d\xi}=0,\label{eq:burgersODE}
\end{equation}
which has the solution (when $\alpha \neq 1$)
\begin{equation}
    -\xi = F + KF^{\frac{\alpha}{\alpha -1}}
\end{equation}
for some $K>0$ (negative $K$ would give a multi-valued function). When $\alpha = 1$, it then follows from (\ref{eq:burgersODE}) that $F = - \xi$, which cannot be matched onto the outer region away from the shock \cite{Eggers09}. In order for the solution to be smooth, $\alpha/(\alpha-1)$ must be an odd positive integer greater than $1$, which gives a discrete set of possible $\alpha$. It can then be shown that $\alpha/(\alpha-1)=3$ is the only stable solution \cite{Eggers09} and hence $\alpha = 3/2$. Then, unravelling definitions,
\begin{equation}
    u (x,t)= u_* + (t_*-t)^{1/2}F\left(\frac{x-x_*+u_*(t_*-t)}{(t_*-t)^{3/2}}\right),
\end{equation}
is the local ($t\rightarrow t_*^-$) similarity solution of the shock formation of the Burgers' equation, where $F = F(\xi)$ is given by
\begin{equation}
    -\xi = F+KF^3.
\end{equation}

In summary, the intuition for the derivation is as follows. First, the Burgers' equation is an advection equation (\ref{eq:advection}) to leading order. Then, the translation to local variables gives an equation local to the shock formation (\ref{eq:burgers_shockframe}) that allows for a self-similar ansatz (\ref{eq:ansatz_burgers}). As a side note, the self-similar solution of the Burgers' equation is of the \textit{second kind}, since the similarity exponent $\alpha$ is only deduced once the equations have been solved.   

\section{Main result}
\label{subsec:main}

The main result of this paper, concerning solutions to (\ref{eq:hyperbolic_PDE}), is stated as follows. Consider a finite-time shock that forms at $(x,t) = (x_*,t_*)$ with $\boldsymbol{f}(x_*,t_*)=\boldsymbol{f}_*$ from smooth initial conditions for a first-order strictly hyperbolic system of PDEs  (\ref{eq:hyperbolic_PDE}). Assume that $\boldsymbol{M}(\boldsymbol{f})$ is smooth at $\boldsymbol{f}=\boldsymbol{f}_*$. Then, the shock formation is locally self-similar as $t \rightarrow t_*^-$ with the solution given by
    \begin{equation}
        \boldsymbol{f}(x,t) = \boldsymbol{f}_* + (t_*-t)^{1/2}F\left(\frac{x-x_*-\lambda(t_*-t)}{c (t_*-t)^{3/2}}\right)\boldsymbol{e},\label{eq:sol_intro}
    \end{equation}
for some eigenvector $\boldsymbol{e}$, with eigenvalue $\lambda$, of $\boldsymbol{M}(\boldsymbol{f}_*)$. Furthermore, the function $F(\xi)$ satisfies $-\xi = F+KF^3$ for some constant $K>0$ and $c$ is a constant that may be derived analytically from $\boldsymbol{M}$ and its derivative (see (\ref{eq:c})).

In other words, the generic shock formations of a first-order system of strictly hyperbolic PDEs (one-dimensional in space) are like that of Burgers' equation and the local similarity solution is  universal, with the shock occurring along an eigenvector $\boldsymbol{e}$ (i.e., $\boldsymbol{f}-\boldsymbol{f}_* \parallel \boldsymbol{e}$). The only unknown is the constant $K>0$, as $(x_*,t_*,\boldsymbol{f}_*)$ are taken as knowns. In general, $(x_*,t_*,\boldsymbol{f}_*)$ are only known numerically and are expected to depend on the initial data. 

Finally it should be noted that there will be higher-order corrections to (\ref{eq:sol_intro}). If the $i$th component of $\boldsymbol{e}$ vanishes, higher-order terms will be important in analyzing the behavior of the $i$th component of $\boldsymbol{f}$. 

\section{Derivation}
\label{sec:deriv}

The derivation of the similarity solution is given in this section. In \cref{subsec:leading_exp}, the leading-order expansion of (\ref{eq:hyperbolic_PDE}) is analysed, which gives an advection equation. Then, the higher-order expansion outlined in \cref{subsec:higher_exp} gives the shock formation. Finally, by imposing that the local similarity solution may be matched to the solution away from the shock, some constants of the solution are evaluated, as described in \cref{subsec:matching}. Thus, the structure of the derivation is closely related to the derivation for the Burgers' equation described in \cref{subsec:background}.

\subsection{Leading-order expansion}\label{subsec:leading_exp}

Let a finite-time shock form at $(x_*,t_*)$ at $\boldsymbol{f}(x_*,t_*)=\boldsymbol{f}_*$. Define variables local to the shock:
\begin{equation}
    \left(x',\tau, \boldsymbol{f}'\right):=\left(x-x_*,t_*-t, \boldsymbol{f}-\boldsymbol{f}_*\right),\label{eq:defn_variables}
\end{equation}
where the sign of $\tau$ is again chosen so that $\tau$ remains positive as $t \rightarrow t_*^-$. Assume that the matrix function $\boldsymbol{M}(\boldsymbol{f})$ in (\ref{eq:hyperbolic_PDE})  is smooth about $\boldsymbol{f} = \boldsymbol{f}_*$. We consider a strictly hyperbolic system of PDEs. Then, expanding to leading (linear) order, (\ref{eq:hyperbolic_PDE}) gives
\begin{equation}
    \frac{\partial f_i'}{\partial \tau}+\sum_{j=1}^N M_{ij}|_{\boldsymbol{f}=\boldsymbol{f}_*}\frac{\partial f_j'}{\partial x'}=0 \quad \forall i = 1,\cdots, N. \label{eq:leading_exp}
\end{equation}
Equation (\ref{eq:leading_exp}) is a linear advection equation. Then, letting $\boldsymbol{e}^{(1)},\cdots, \boldsymbol{e}^{(N)}$ be the basis of eigenvectors of $\boldsymbol{M}|_{\boldsymbol{f}=\boldsymbol{f}_*}$ with corresponding (distinct) eigenvalues $\lambda^{(1)},\cdots, \lambda^{(N)}$, the solution is given as 
\begin{equation}
    \boldsymbol{f}'(x',\tau) = \sum_{a=1}^N g^{(a)}\left(x'-\lambda^{(a)}\tau\right)\boldsymbol{e}^{(a)} \label{eq:leading_order_soln}
\end{equation}
for some functions $g^{(1)},\cdots, g^{(N)}$. 

As the characteristics are constant, the solution to the linear advection equation (\ref{eq:leading_order_soln}) cannot form a singularity from smooth conditions (see Appendix \ref{app:smooth_advection}). Thus, a shock does not form in the leading-order expansion and the analysis must be continued to higher order to describe analytically the solution as the shock forms. 

\subsection{Higher-order expansion}\label{subsec:higher_exp}

The expansion accounting for the next higher-order term in (\ref{eq:hyperbolic_PDE}) is given by (again, smoothness of $\boldsymbol{M}$ is assumed)
\begin{equation}
    \frac{\partial f_i'}{\partial \tau}+\sum_{j=1}^N M_{ij}|_{\boldsymbol{f}=\boldsymbol{f}_*}\frac{\partial f_j'}{\partial x'}+\sum_{j=1}^N \sum_{k=1}^N M_{ij, k}|_{\boldsymbol{f}=\boldsymbol{f}_*}f_k'\frac{\partial f_j'}{\partial x'}=0 \quad \forall i = 1,\cdots, N, \label{eq:nord_exp}
\end{equation}
where
\begin{equation}
    M_{ij, k}|_{\boldsymbol{f}=\boldsymbol{f}_*}:=\left.\frac{\partial M_{ij}}{\partial f_k}\right|_{\boldsymbol{f} = \boldsymbol{f}_*}.
\end{equation}

From the leading-order equations (\ref{eq:leading_order_soln}), it is expected that a shock forms with the velocity given by some eigenvalue $\lambda$ (denote the corresponding eigenvector by $\boldsymbol{e}$). It is natural to then consider the reference frame of the shock by defining the variable    
\begin{equation}
    \overline{x} := x' - \lambda \tau \label{eq:shock_coords}.
\end{equation}
Then, in $(\overline{x},\tau)$ coordinates, (\ref{eq:nord_exp}) gives
\begin{equation}
    \frac{\partial f_i'}{\partial \tau}+\sum_{j=1}^N\left(M_{ij}|_{\boldsymbol{f}=\boldsymbol{f}_*}-\lambda \delta_{ij}\right)\frac{\partial f_j'}{\partial \overline{x}}+\sum_{j=1}^N \sum_{k=1}^N M_{ij, k}|_{\boldsymbol{f}=\boldsymbol{f}_*}f_k'\frac{\partial f_j'}{\partial \overline{x}}=0 \quad \forall i = 1,\cdots, N, \label{eq:next_order_shock_frame}
\end{equation}
where $\delta_{ij}$ is the Kronecker delta. As for the Burgers' equation, the scalings for the terms in (\ref{eq:next_order_shock_frame}) are as follows. Let $\overline{x} = \textit{O}(\tau^{\alpha})$ as $\tau \rightarrow 0^+$. Let the maximum magnitude component as $\tau \rightarrow 0^+$ of $\boldsymbol{f}'$ be $\textit{O}(\tau^{\beta})$. The maximum magnitude component of $\boldsymbol{f}' \partial \boldsymbol{f}'/\partial \overline{x}$ is then $\textit{O}(\tau^{2\beta-\alpha})$. Since  $\boldsymbol{f}' = 0$ at $\tau = 0$, it follows that $\beta >0$ and hence higher order correction terms to the expansion (\ref{eq:next_order_shock_frame}) such as the components of $\boldsymbol{f}' \boldsymbol{f}' \partial \boldsymbol{f}'/\partial \overline{x}$, will be subdominant. Since the maximum component of $\partial \boldsymbol{f}'/\partial \tau$ is $\textit{O}(\tau^{\beta-1})$, dominant balance gives $\beta = \alpha - 1$. Note that the scaling $\boldsymbol{f}' = \textit{O}(\tau^{\alpha-1})$ is for the component with maximum magnitude and thus it is possible for some components to be zero at $\textit{O}(\tau^{\alpha-1})$ order, where a simple example is shock formation of the system with two dependent variables $u,v$ with $u$ satisfying the Burgers' equation (\ref{eq:burger}) and $v$ satisfying the equation $\partial v/\partial t = 0$.
The scalings identified implies that the leading-order term in (\ref{eq:next_order_shock_frame}) is order $\textit{O}(\tau^{-1})$ and satisfies
$\sum_{j=1}^N\left(M_{ij}|_{\boldsymbol{f}=\boldsymbol{f}_*}-\lambda \delta_{ij}\right)\partial f_j'/\partial \overline{x}=0 ~\forall i = 1,\cdots, N$. In other words, since the eigenvalues are distinct, $\partial \boldsymbol{f}'/\partial \overline{x}$ is parallel to $\boldsymbol{e}$. 

When the system of PDEs is not strictly hyperbolic (i.e., repeated eigenvalues), it is possible that $\partial \boldsymbol{f}'/\partial \overline{x}$ is a linear combination of linearly independent eigenvectors with the same eigenvalue $\lambda$. This is the only step in the derivation where the uniqueness of eigenvalues is used, and hence the method is also applicable to hyperbolic systems with repeated eigenvalues if it is true that $\partial \boldsymbol{f}'/\partial \overline{x}$ is parallel to some eigenvector $\boldsymbol{e}$.

It then follows that the expansion for $\boldsymbol{f}'$ as $\tau \rightarrow 0^+$ is given by
\begin{equation}
    \boldsymbol{f}'(\overline{x},\tau) = g(\overline{x}, \tau) \boldsymbol{e} + \boldsymbol{q} \tau^{\alpha-1} + \boldsymbol{h}(\overline{x}, \tau), \label{eq:expansion}
\end{equation}
where $g = \textit{O}(\tau^{\alpha-1})$, $\boldsymbol{h} = \textit{O}(\tau^{2\alpha-2})$ is a higher-order correction term, and $\boldsymbol{q}$ is some constant. The scaling for $\boldsymbol{h}$ can be seen by balancing $(\partial g/\partial \tau)\boldsymbol{e}$ with $\partial \boldsymbol{h}/\partial \overline{x}$. Since $\boldsymbol{M}|_{\boldsymbol{f}=\boldsymbol{f}_*}\cdot \boldsymbol{e} - \lambda \boldsymbol{e} = \boldsymbol{0}$, then substituting the expansion (\ref{eq:expansion}) into (\ref{eq:next_order_shock_frame}) gives 
\begin{multline}
    \frac{\partial}{\partial \tau}(e_i g(\overline{x},\tau) + q_i \tau^{\alpha-1})+\sum_{j=1}^N \left(M_{ij}|_{\boldsymbol{f}=\boldsymbol{f}_*}-\lambda \delta_{ij}\right)  \frac{\partial h_j}{\partial \overline{x}}\\
    +\sum_{j=1}^N \sum_{k=1}^N M_{ij,k}|_{\boldsymbol{f}=\boldsymbol{f}_*}(e_k g(\overline{x},\tau) + q_k \tau^{\alpha-1})\frac{\partial g}{\partial \overline{x}}e_j =0 \quad \forall i = 1,\cdots, N.\label{eq:substituted_expansion}
\end{multline}

Before proceeding further, it should be observed that (\ref{eq:substituted_expansion}) is a system of $N$ equations for $N+1$ functions $g$ and $h_1, \cdots, h_N$. However, there is a degeneracy in the coefficient of the $\partial h_j /\partial \overline{x}$ term, since $\det(\boldsymbol{M}|_{\boldsymbol{f}=\boldsymbol{f}_*}-\lambda \boldsymbol{I})=0$ as $\lambda$ is an eigenvalue of $\boldsymbol{M}|_{\boldsymbol{f}=\boldsymbol{f}_*}$. Then, (\ref{eq:substituted_expansion}) can be simplified by taking a suitable linear combination of the $N$ equations. Indeed, there exists a left eigenvector $\boldsymbol{e}^L$ of $\boldsymbol{M}|_{\boldsymbol{f}=\boldsymbol{f}_*}$ with eigenvalue $\lambda$ satisfying (see \cref{app:left_evec})
\begin{equation}
    \boldsymbol{e}^L \cdot \boldsymbol{M}|_{\boldsymbol{f}=\boldsymbol{f}_*} -\lambda \boldsymbol{e}^L=\boldsymbol{0} \text{ and } \boldsymbol{e}^L \cdot \boldsymbol{e}\neq 0. \label{eq:left_evec_condition}
\end{equation}
Upon taking the inner product of (\ref{eq:substituted_expansion}) with $\boldsymbol{e}^L$, it then follows that
\begin{multline}
    \sum_{i=1}^N\frac{\partial}{\partial \tau}\left (e_ie^L_i g(\overline{x},\tau) + q_ie^L_i \tau^{\alpha-1}\right )\\+ \sum_{i=1}^N\sum_{j=1}^N \sum_{k=1}^N M_{ij,k}|_{\boldsymbol{f}=\boldsymbol{f}_*}\left (e^L_ie_k g(\overline{x},\tau) + e^L_iq_k \tau^{\alpha-1}\right )\frac{\partial g}{\partial \overline{x}}e_j =0\label{eq:subsituted_expansion_simplified}.
\end{multline}

Finally, there is a degree of freedom in the expansion (\ref{eq:expansion}) for the choice of $(g,\boldsymbol{q})$. Then, without loss of generality, we consider $\boldsymbol{q} \cdot \boldsymbol{e}^L=0$ since $(g,\boldsymbol{q})$ can be redefined through the projection $(g + (\boldsymbol{q}\cdot \boldsymbol{e}^L)(\boldsymbol{e}\cdot \boldsymbol{e}^L)^{-1},\boldsymbol{q}-(\boldsymbol{q}\cdot \boldsymbol{e}^L)(\boldsymbol{e}\cdot \boldsymbol{e}^L)^{-1} \boldsymbol{e})$. Thus, (\ref{eq:subsituted_expansion_simplified}) gives that
\begin{equation}
    \frac{\partial g}{\partial \tau}+c_1 \tau^{\alpha-1} \frac{\partial g}{\partial \overline{x}}-c_2 g \frac{\partial g}{\partial \overline{x}}=0,\label{eq:before_scaling_burgers}
\end{equation}
where there are two constants
\begin{subeqnarray}
    c_1 &=& \frac{\sum_{i=1}^N\sum_{j=1}^N \sum_{k=1}^N M_{ij,k}|_{\boldsymbol{f}=\boldsymbol{f}_*}e^L_i e_j q_k}{\sum_{i=1}^Ne^L_ie_i},\\
    c_2 &=&-\frac{\sum_{i=1}^N\sum_{j=1}^N \sum_{k=1}^N M_{ij,k}|_{\boldsymbol{f}=\boldsymbol{f}_*}e^L_i e_j e_k}{\sum_{i=1}^Ne^L_ie_i}.
\end{subeqnarray}

In this work, we assume that $c_2 \neq 0$. The degenerate case $c_2 =0$ would require the consideration of higher-order terms, which is not considered. Rescaling via $\overline{x}_s := c_2^{-1}\left(\overline{x} - c_1 \alpha^{-1}\tau^{\alpha}\right)$, the Burgers' equation (\ref{eq:burgers_shockframe}) is recovered: 
\begin{equation}
    \frac{\partial g}{\partial \tau}- g\frac{\partial g}{\partial \overline{x}_s}=0.
\end{equation}
Thus, as seen in \cref{subsec:background}, $\alpha = 3/2$ and $g = \tau^{1/2}F(\xi)$ where $\xi :=\overline{x}_s \tau^{-3/2}$ and
\begin{equation}
    \xi = -F - KF^3 \label{eq:alg_eq}
\end{equation}
for some $K>0$. We have therefore found the similarity solution in the neighborhood of the shock.

\subsection{Matching}\label{subsec:matching}

Having found the similarity solution in \cref{subsec:higher_exp}, we may go one step further and show that the vector $\boldsymbol{q} = \boldsymbol{0}$ by requiring that the similarity solution local to the shock must match to the solution away from the shock. Explicitly, such a condition is given by
\begin{equation}
    \lim_{\tau \rightarrow 0^+}\boldsymbol{f}'(\Delta, \tau) = \boldsymbol{f}'(\Delta,0),
\end{equation}
for any fixed $\overline{x} = \Delta$ \cite{Eggers09, EggersFontelos}.

Since $g(\Delta, \tau) = \tau^{1/2}F(c_2^{-1}\left(\Delta - (2c_1/3)\tau^{3/2}\right) \tau^{-3/2})$ and $F(\xi) \sim -\xi^{1/3}$ for $|\xi|$ large (from (\ref{eq:alg_eq})), for any fixed $\overline{x} = \Delta$, we have that 
\begin{equation}
    \lim_{\tau \rightarrow 0^+}g(\Delta, \tau) \text{ is bounded.} \label{eq:g_bounded}
\end{equation}
Thus, from the expansion (\ref{eq:expansion}), we have that
\begin{equation}
    \lim_{\tau \rightarrow 0^+}\boldsymbol{h}(\Delta, \tau) \text{ is bounded}.\label{eq:h_bounded}
\end{equation}

We may use (\ref{eq:g_bounded}) and (\ref{eq:h_bounded}) to deduce the value of $\boldsymbol{q}$. Indeed, substituting (\ref{eq:before_scaling_burgers}) into (\ref{eq:substituted_expansion}) gives that
\begin{multline}
    -c_1e_i \tau^{\frac{1}{2}}\frac{\partial g}{\partial \overline{x}}-c_2 e_i g \frac{\partial g}{\partial \overline{x}}+\sum_{j=1}^N\left(M_{ij}|_{\boldsymbol{f}=\boldsymbol{f}_*}-\lambda \delta_{ij}\right)  \frac{\partial h_j}{\partial \overline{x}} \\
    + \sum_{j=1}^N \sum_{k=1}^NM_{ij,k}|_{\boldsymbol{f}=\boldsymbol{f}_*}(e_k g(\overline{x},\tau) + q_k \tau^{\frac{1}{2}})\frac{\partial g}{\partial \overline{x}}e_j =-\frac{1}{2}q_i\tau^{-\frac{1}{2}} \quad \forall i = 1,\cdots, N,
\end{multline}
and integration with respect to $\overline{x}$ gives
\begin{multline}
    G_i(\tau)-c_1e_i \tau^{\frac{1}{2}}g-c_2 e_i \frac{g^2}{2}+\sum_{j=1}^N\left(M_{ij}|_{\boldsymbol{f}=\boldsymbol{f}_*}-\lambda \delta_{ij}\right)  h_j \\
    + \sum_{j=1}^N \sum_{k=1}^N M_{ij,k}|_{\boldsymbol{f}=\boldsymbol{f}_*}\left(\frac{1}{2}e_k g(\overline{x},\tau) + q_k \tau^{\frac{1}{2}}\right)ge_j =-\frac{1}{2}q_i\overline{x}\tau^{-\frac{1}{2}} \quad \forall i = 1,\cdots, N,\label{eq:integrated_qi_eq}
\end{multline}
for some function $G_i(\tau)$. Now, we have from (\ref{eq:g_bounded}) and (\ref{eq:h_bounded}) that the left-hand side of (\ref{eq:integrated_qi_eq}) is bounded as $\tau \rightarrow 0^+$ for any fixed $\overline{x} = \Delta$ and hence $q_i = 0$ for all $i=1,\cdots, N$. 

Then, the similarity solution is given by
\begin{equation}
    \boldsymbol{f}(x,t) = \boldsymbol{f}_* + (t_*-t)^{1/2}F\left(\frac{x-x_*-\lambda(t_*-t)}{c (t_*-t)^{3/2}}\right)\boldsymbol{e},\label{eq:soln}
\end{equation}
where
\begin{equation}
    c = -\frac{\sum_{i=1}^N\sum_{j=1}^N \sum_{k=1}^N M_{ij,k}|_{\boldsymbol{f}=\boldsymbol{f}_*}e^L_i e_j e_k}{\sum_{i=1}^Ne^L_ie_i}. \label{eq:c}
\end{equation}
and $F(\xi)$ satisfies (\ref{eq:alg_eq}). From the expansion (\ref{eq:expansion}), corrections to (\ref{eq:soln}) are $\textit{O}(t_*-t)$.

Throughout the derivation, we did not specify a magnitude for $\boldsymbol{e}$. Indeed, the solution (\ref{eq:alg_eq}, \ref{eq:soln}, \ref{eq:c}) is exactly invariant under the transformation $(\boldsymbol{e},c, \xi, F, K)$ to $(\mu \boldsymbol{e}, \mu c, \mu^{-1}\xi, \mu^{-1}F, \mu^2K)$ for constant $\mu \neq 0$. 

\section{Corollaries}\label{sec:corollaries}

Due to the relatively simple algebraic equation (\ref{eq:alg_eq}) that the similarity function $F$ satisfies, we can write down some practical expressions, which may be useful when using the presented method in physical problems. Since there are two exponents in the similarity solution ($\boldsymbol{f}' \sim (t_*-t)^{1/2}, \overline{x} \sim (t_*-t)^{3/2}$), we derive two power-law expressions. The expressions are used in the numerical cross-check of the solution, as shown in \cref{sec:verif}.

First, we write down how the first derivatives become infinite as the shock singularity is approached. That is, as $t \rightarrow t_*^-$,
\begin{equation}
    \max \left|\frac{\partial f_i}{\partial x}\right| = \frac{|e_i c^{-1}|}{(t_*-t)} \max \left|\frac{dF}{d\xi}\right| = \frac{|e_i c^{-1}|}{(t_*-t)} \quad \forall i=1,\cdots, N\label{eq:first_deriv},
\end{equation}
where the second equality follows since $dF/d\xi=-(1+3KF^2)^{-1}$. The expression (\ref{eq:first_deriv}) does not contain any unknowns (the values $x_*,t_*, \boldsymbol{f}_*$ are taken as known), and is consequently useful for cross-checking with numerical solutions. 

Similarly, we consider how the second derivatives become infinite as the shock is approached. As $t \rightarrow t_*^-$,
\begin{equation}
    \max \left|\frac{\partial^2 f_i}{\partial x^2}\right| = \frac{|e_ic^{-2}|}{(t_*-t)^{5/2}}\frac{25\sqrt{15}}{108}K^{1/2} \quad \forall i=1,\cdots, N,\label{eq:second_deriv}
\end{equation}
which follows from evaluating $d^2F/d\xi^2$ using the algebraic equation (\ref{eq:alg_eq}). The expression (\ref{eq:second_deriv}) contains the unknown constant $K$. However the constant $K$ deduced from fitting (\ref{eq:second_deriv}) can then be used in checking the self-similar collapse of the $\boldsymbol{f}$ profiles, which therefore means that $K$ can also be independently checked (see \cref{sec:verif}).

\section{Example}
\label{sec:verif}

Below, we present an example for clarity and to independently check the theoretical predictions with numerical solutions. The method presented in this paper is applied to the (nondimensional) shallow water equations in one-dimensional space, given by
\begin{subeqnarray}
    \frac{\partial u}{\partial t}&=&-u\frac{\partial u}{\partial x}-\frac{\partial \eta}{\partial x},\\
    \frac{\partial \eta}{\partial t}&=&  - \eta \frac{\partial u}{\partial x}-u \frac{\partial \eta}{\partial x},\label{eq:swe_nondim}
\end{subeqnarray}
for velocity $u(x,t)$ and water height $\eta(x,t)$. We apply the method presented in this paper for $\boldsymbol{f} = (u,\eta)$. Equation (\ref{eq:swe_nondim}) is indeed of form (\ref{eq:hyperbolic_PDE}), with
\begin{equation}
    \boldsymbol{M}(\boldsymbol{f}) =\begin{pmatrix}
        -u& -1\\
        -\eta&-u
    \end{pmatrix}.
\end{equation}

The shallow water equations (\ref{eq:swe_nondim}) are solved numerically using centered finite differences and time-stepped with the  Crank-Nicolson approach. A fixed non-uniform grid is used with an adaptive time step. A uniform grid and time step also works well, but the excess computational cost means that the dynamics is harder to resolve close to the singularity at time $t_*$. The code and data are available at \cite{eshima_data_2026}. The numerical results presented in this manuscript have been checked for numerical convergence.  

\begin{figure}[htbp]
  \centering
  \includegraphics[width=\textwidth]{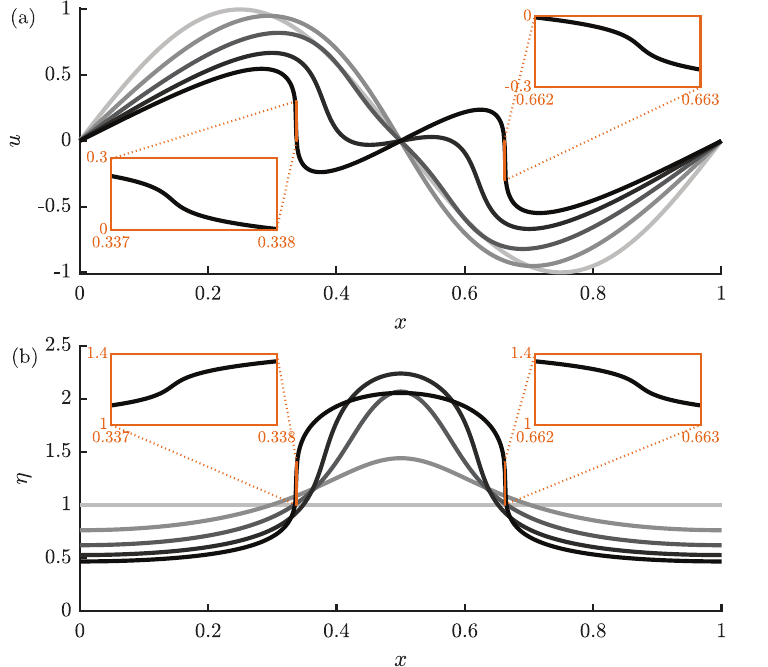}
  \caption{Example evolution of the shallow water equations (\ref{eq:swe_nondim}), taking periodic boundary conditions for the domain $x \in [0,1]$ with initial condition $(u(x,0),\eta(x,0))=(\sin(2\pi x), 1)$. (a) Velocity $u(x,t)$, where the evolution from gray to black shows evolution in time of ($t = 0,0.05, 0.1, 0.15, 0.195$). The shock singularity occurs at $t_* \approx 0.196$ (thus, the final timestep shown is just before the shock forms). Inset shows magnified views of the shock region at $t = 0.195$. (b) Analogous to (a) for the water height $h(x,t)$.}
  \label{fig:sample_evol}
\end{figure}

As the example, we consider the shock formation due to the initial conditions $(u(x,0),\eta(x,0))=(\sin(2\pi x), 1)$ on the domain $x \in [0,1]$, subject to periodic boundary conditions. A shock forms at $t_* \approx 0.196$. \Cref{fig:sample_evol} presents the evolution of $(u,\eta)$ over time. The color (grey to black) shows the evolution from the initial state $t=0$ to the state close to shock formation $t = 0.195$, where the shock singularity (the derivatives become infinite) occurs at $t_* \approx 0.196$. The inset shows a magnified view of the shock region at $t = 0.195$. Due to the initial conditions and symmetry of the equations, we see that the solution is mirror symmetric about $x = \frac{1}{2}$ ($u\leftrightarrow-u$ and $\eta \leftrightarrow \eta$). Then, two shocks form at the same time (one at location $x_* \approx 0.663$ and the other at location $x_* \approx 0.337$). 

Let the shock occur at $\boldsymbol{f}_* = (u_*,\eta_*)$. Then, following the notation as presented in \cref{sec:deriv}:
\begin{equation}
    \boldsymbol{M}|_{\boldsymbol{f}=\boldsymbol{f}_*} =\begin{pmatrix}
        -u_*& -1\\
        -\eta_*&-u_*
    \end{pmatrix},
\end{equation}
which has eigenvectors and eigenvalues
\begin{subeqnarray}
    \boldsymbol{e}^{(\pm)} = (1, \mp \sqrt{\eta_*}), ~\lambda^{(\pm)}=-u_* \pm \sqrt{\eta_*}.
\end{subeqnarray}

It can be observed from numerical data that the shock formation around $x_* \approx 0.663$ is along the $\boldsymbol{e}^{(-)}$ eigenvector, since both $u$ and $\eta$ are decreasing close to $x_*$. For the mirror symmetric shock formation around $x_* \approx 0.337$, the shock forms along the $\boldsymbol{e}^{(+)}$ eigenvector and the example therefore shows the two possibilities. 

For the application of the method presented in this paper, we will focus on the local region around $x_* \approx 0.663$ as the other shock is mirror symmetric. Thus we consider shock formation about $\boldsymbol{e} = \boldsymbol{e}^{(-)}$. 

We calculate the various expressions in the solution (\ref{eq:alg_eq}, \ref{eq:soln}, \ref{eq:c}). The left eigenvector of interest is given by
\begin{equation}
    \boldsymbol{e}^L =  (\sqrt{\eta_*}, 1)
\end{equation}
and the components of $M_{ij,k}|_{\boldsymbol{f}=\boldsymbol{f}_*}$ are given by
\begin{align}
M_{11,1}|_{\boldsymbol{f}=\boldsymbol{f}_*} &= -1, &\quad M_{12,1}|_{\boldsymbol{f}=\boldsymbol{f}_*} &= 0, 
   &\quad M_{21,1}|_{\boldsymbol{f}=\boldsymbol{f}_*} &= 0, &\quad M_{22,1}|_{\boldsymbol{f}=\boldsymbol{f}_*} &= -1, \\
M_{11,2}|_{\boldsymbol{f}=\boldsymbol{f}_*} &=  0, &\quad M_{12,2}|_{\boldsymbol{f}=\boldsymbol{f}_*} &= 0,
   &\quad M_{21,2}|_{\boldsymbol{f}=\boldsymbol{f}_*} &= -1, &\quad M_{22,2}|_{\boldsymbol{f}=\boldsymbol{f}_*} &=  0. \nonumber
\end{align}
Then, (\ref{eq:c}) gives
\begin{equation}
    c = -\frac{-\left(\sqrt{\eta_*} + \sqrt{\eta_*} + \sqrt{\eta_*}\right)}{2\sqrt{\eta_*} } = \frac{3}{2}.
\end{equation}
Using (\ref{eq:soln}), we deduce that as $t \rightarrow t_*^-$,
\begin{subeqnarray}
    u(x,t) &=& u_* +(t_*-t)^{\frac{1}{2}}F\left(\frac{x-x_*-\lambda(t_*-t)}{\frac{3}{2}(t_*-t)^{\frac{3}{2}}}\right),\\
    \eta(x,t) &=& \eta_* +\sqrt{\eta_*}(t_*-t)^{\frac{1}{2}}F\left(\frac{x-x_*-\lambda(t_*-t)}{\frac{3}{2}(t_*-t)^{\frac{3}{2}}}\right).\label{eq:simil_soln_example}
\end{subeqnarray}
From (\ref{eq:first_deriv}) we deduce that as $t \rightarrow t_*^-$,
\begin{equation}
    \left(\max \left|\frac{\partial u}{\partial x}\right|,  \max \left|\frac{\partial \eta}{\partial x}\right| \right) = \frac{2}{3(t_*-t)}(1,\sqrt{\eta_*})\label{eq:first_deriv_example}
\end{equation}
and from (\ref{eq:second_deriv}) that as $t \rightarrow t_*^-$,
\begin{equation}
    \left(\max \left|\frac{\partial^2 u}{\partial x^2}\right|, \max \left|\frac{\partial^2 \eta}{\partial x^2}\right| \right)= \frac{1}{(t_*-t)^{5/2}}\frac{25\sqrt{15}}{243}K^{1/2}(1,\sqrt{\eta_*}) \label{eq:second_deriv_example}
\end{equation}

\begin{figure}[htbp]
  \centering
  \includegraphics[width=\textwidth]{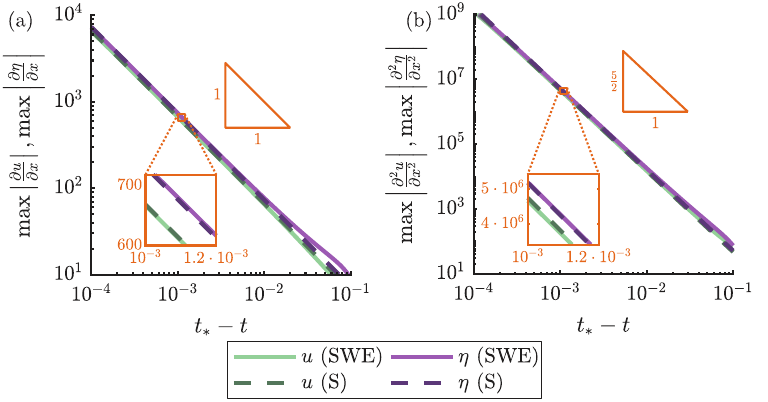}
  \caption{
  Comparison of the theoretical and numerical similarity solution exponents and prefactors. The shallow water equations (\ref{eq:swe_nondim}) are solved numerically, taking periodic boundary conditions for the domain $x \in [0,1]$ with initial condition $(u(x,0),\eta(x,0))=(\sin(2\pi x), 1)$. The numerical solution of the shallow water equations (SWE) (solid curves) is compared against the similarity solution (S) (dashed lines) predicted by the similarity solution (\ref{eq:first_deriv_example}, \ref{eq:second_deriv_example}). (a,b) Log-log plot of $\max |\partial u/\partial x|$, $\max |\partial \eta/\partial x|$ (green, purple) as the shock forms $t \rightarrow t_*^-$, showing a $(t_*-t)^{-1}$ blowup. There are no fitting parameters for the prediction (\ref{eq:first_deriv_example}). (b)  Log-log plot of $\max |\partial^2 u/\partial x^2|$, $\max |\partial^2 \eta/\partial x^2|$ (green, purple) as the shock forms $t \rightarrow t_*^-$, showing a $(t_*-t)^{-5/2}$ blowup, for the best fit $K \approx 0.14$, the only unknown in the prediction (\ref{eq:second_deriv_example}). Insets show a magnified view for $t_*-t \in [10^{-3}, 1.2 \times 10^{-3}]$.}
  \label{fig:verif}
\end{figure}

\Cref{fig:verif} compares the theoretical and numerical predictions for the exponents and prefactors of the similarity solutions for the expressions (\ref{eq:first_deriv_example}, \ref{eq:second_deriv_example}). It can be seen that the numerical solution (solid) is in excellent agreement with the similarity solutions (dashed) given in (\ref{eq:first_deriv_example}, \ref{eq:second_deriv_example}). The analytical prediction (\ref{eq:first_deriv_example}) used for \cref{fig:verif}(a) does not have any unknowns (i.e., no fitting parameters), whereas the analytical prediction (\ref{eq:second_deriv_example}) for \cref{fig:verif}(b) contains the single unknown $K$ (i.e., one fitting parameter). For the numerical solution, the singularity values $(x_*,t_*,u_*,\eta_*)$ can be obtained by extrapolation of the values of $(x,t,u,\eta)$ as the shock is formed and are therefore considered known. The figure presents the numerical solutions down to $t_*-t = 10^{-4}$, which is the threshold for which our numerical solution is converged. Improved numerical accuracy would allow the results to be converged for smaller values of $t_*-t$. 

\begin{figure}[htbp]
  \centering
  \includegraphics[width=\textwidth]{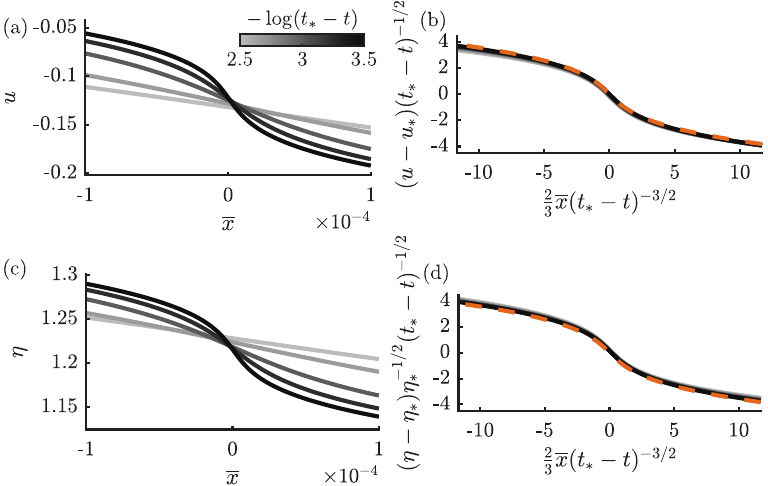}
  \caption{
  Comparison of the theoretical and numerical similarity solution profiles. The shallow water equations (\ref{eq:swe_nondim}) are solved numerically, taking periodic boundary conditions for the domain $x \in [0,1]$ with initial condition $(u(x,0),\eta(x,0))=(\sin(2\pi x), 1)$. The numerical solutions (solid curves) are colored gray to black as $t\rightarrow t_*^-$. The reference frame $\overline{x} = x-x_*-\lambda (t_*-t)$ of the shock is considered. The orange dashed curves show the similarity solution prediction $-X = Y+KY^3$ where $X$ is the horizontal axis and $Y$ is the vertical axis for $K \approx 0.14$ as predicted in \cref{fig:verif}(b). (a) $u$ against $\overline{x}$. (b) $u$ against $\overline{x}$, suitably rescaled according to the similarity solution (\ref{eq:simil_soln_example}). (c,d) Analogous to (a,b) for the water height $\eta$.}
  \label{fig:collapse}
\end{figure}

We may cross-check the value of $K$ predicted from $\partial^2 u/\partial x^2$ and $\partial^2 \eta/\partial x^2$ in \cref{fig:verif}(b) by considering the self-similar profiles. \Cref{fig:collapse} compares the theoretical and numerical profiles of $u$ and $\eta$. The numerical solutions are given by the solid curves (colored gray to black as $t \rightarrow t_*^-$) and the dashed curves (orange) show the theoretical prediction, using the value $K$ predicted from $\partial^2 u/\partial x^2$ and $\partial^2 \eta/\partial x^2$. Once again, there is excellent agreement.

\section{Conclusions}
\label{sec:conclusions}

In this paper, we have given a method to derive the self-similar solutions of one-dimensional strictly hyperbolic PDEs, by generalizing the solution method for the Burgers' equation. Given the universality and simplicity of the method, it should serve as a valuable tool for future analysis of hyperbolic PDEs. 

As direct extension, it would be interesting to consider the case of repeated eigenvalues of $\boldsymbol{M}$. In addition, it would be insightful to analyze self-similarity for multiple dimensions in space, such as in \cite{Eggers09} for the compressible gas equations.

\appendix

\section{Smoothness of linear advection equation} \label{app:smooth_advection}

For the solution (\ref{eq:leading_order_soln}), suppose that $\boldsymbol{f}'$ is smooth for $\tau >0$. Recall that the eigenvalues $\lambda^{(a)}$ are distinct. Also, note the following identity:
\begin{multline}
\left(\boldsymbol{M}|_{\boldsymbol{f}=\boldsymbol{f}_*} - \lambda^{(2)} \boldsymbol{I}\right)\cdot \left(\boldsymbol{M}|_{\boldsymbol{f}=\boldsymbol{f}_*} - \lambda^{(3)} \boldsymbol{I}\right)\cdots \left(\boldsymbol{M}|_{\boldsymbol{f}=\boldsymbol{f}_*} - \lambda^{(N)} \boldsymbol{I}\right)\cdot \boldsymbol{f}' \\
= \left(\prod_{a = 2}^{N}\left(\lambda^{(1)}-\lambda^{(a)}\right)\right) g^{(1)}\left(x'-\lambda^{(1)}\tau\right)\boldsymbol{e}^{(1)}.
\end{multline}
Then, it follows that $g^{(1)}(x'-\lambda^{(1)}\tau)$ is also smooth for $\tau >0$. By considering the coordinate $x'':= x'-\lambda^{(1)}\tau$, it then follows that  $g^{(1)}(x'-\lambda^{(1)}\tau)$ is smooth at $\tau = 0$. By the same argument, all $g^{(a)}$ are smooth at $\tau = 0$ and hence $\boldsymbol{f}'$ also. Thus, a singularity cannot form for the solution of the linear advection equation (\ref{eq:leading_order_soln}) from smooth conditions. 

\section{Left eigenvectors}\label{app:left_evec} 

The matrix $\boldsymbol{M}|_{\boldsymbol{f}=\boldsymbol{f}_*}$ is diagonalisable with $\boldsymbol{M}|_{\boldsymbol{f}=\boldsymbol{f}_*} = \boldsymbol{P} \boldsymbol{\Lambda} \boldsymbol{P}^{-1}$ for diagonal matrix $\boldsymbol{\Lambda}= \text{diag}(\lambda^{(1)},\lambda^{(2)},\dots, \lambda^{(N)})$ and $\boldsymbol{P}$ is the matrix whose $i$th column is given by $\boldsymbol{e}^{(i)}$.
Let $\boldsymbol{w}^{(i)}$ be the $i$th row of $\boldsymbol{P}^{-1}$, which are the left eigenvalues of $\boldsymbol{M}|_{\boldsymbol{f}=\boldsymbol{f}_*}$.
Then, since $\boldsymbol{P}^{-1}\boldsymbol{P} =\boldsymbol{I}$, we have that for any $i$, $\boldsymbol{w}^{(i)} \cdot \boldsymbol{e}^{(i)} = 1$. Then, $\boldsymbol{w}^{(i)}$ and any non-zero vector parallel to $\boldsymbol{w}^{(i)}$, has the required properties given in (\ref{eq:left_evec_condition}).

\section*{Acknowledgments} We thank Peter Constantin, Jens Eggers, and Marco Fontelos for helpful discussions.
\bibliographystyle{siamplain}
\bibliography{references}
\end{document}